\input amssym
\font\sevenbf=ptmb at 7.6pt

\font\sevensl=ptmro at 7.6pt
\font\sevenrm=ptmr at 7.6pt
\rm

\font\tt=cmtt10

\textfont"F=\teni
 \scriptfont"F=\seveni
 \scriptscriptfont"F=\fivei

 \mathchardef\flat="2F5B
 \mathchardef\natural="2F5C
 \mathchardef\sharp="2F5D
 \mathchardef\star="2F3F
 \mathchardef\triangleleft="2F2F

\font\mathlet=eurm10
 \font\mathsublet=eurm7
 \font\mathsubsublet=eurm5
\textfont1=\mathlet
 \scriptfont1=\mathsublet
 \scriptscriptfont1=\mathsubsublet
 \skewchar\mathlet='177 \skewchar\mathsublet='177 \skewchar\mathsubsublet='177
\mathcode`0="7130
\mathcode`1="7131
\mathcode`2="7132
\mathcode`3="7133
\mathcode`4="7134
\mathcode`5="7135
\mathcode`6="7136
\mathcode`7="7137
\mathcode`8="7138
\mathcode`9="7139
\mathchardef\Gamma="7100
\mathchardef\Delta="7101
\mathchardef\Theta="7102
\mathchardef\Lambda="7103
\mathchardef\Xi="7104
\mathchardef\Pi="7105
\mathchardef\Sigma="7106
\mathchardef\Upsilon="7107
\mathchardef\Phi="7108
\mathchardef\Psi="7109
\ifx\eplain\undefined \input eplain \fi

\newcount\subsecno\subsecno=0
\newcount\thmno\thmno=0
\newcount\secno\secno=0

\outer\def\section #1\par{
\vskip 12pt plus.1\vsize\penalty-350
\vskip 0pt plus-.1\vsize\bigskip\vskip\parskip
\global\advance\secno by 1\global\thmno=0\global\subsecno=0
\def\label##1{\definexref{##1}{\the\secno}{Part}}
{\teneufm\the\secno.\ }{\bf #1}\par
\smallskip
}

\outer\def\subsection #1\par{
\vskip 0pt plus.1\vsize\penalty-250
\vskip 0pt plus-.1\vsize\bigskip\vskip\parskip
\global\advance\subsecno by 1\global\thmno=0
\def\label##1{\definexref{##1}{\the\secno.\the\subsecno}{Section}}
{\teneufm\the\secno.\the\subsecno.\ }{\bf #1}\par
\smallskip
}

\def\tag{\the\secno.\the\subsecno.\the\thmno}

\everydisplay{\def\label#1{
\global\advance\thmno by 1
\definexref{#1}{(\tag)}{Equation}
\leqno(\tag)}}

\def\myproclaim#1#2#3#4{
\global\advance\thmno by 1\def\label##1{\definexref{##1}{\tag}{#1}}
\medbreak\noindent{#3#1 \tag.}\enspace{#4#2}\par
\ifdim\lastskip<\medskipamount \removelastskip\penalty55\medskip\fi}
\outer\def\thm #1. #2\par{\myproclaim{#1}{#2}{\bf}{\sl}}
\outer\def\rmk #1. #2\par{\myproclaim{#1}{#2}{\it}{}}
\mathchardef\Omega="710A
\magnification1200
\nopagenumbers

\font\sevenss=cmss7
\headline{\everymath{\scriptstyle}\sevenss \ifnum\pageno>1
\ifodd\pageno\title\ \hrulefill\ \the\pageno
\else\the\pageno\ \hrulefill\ S.~R.~Gal\fi\fi}

\def\beginabstract{\bigskip\midinsert
\narrower\narrower\noindent\everymath{\scriptstyle}%
\sevenbf Abstract: \sevenrm}

\def\endabstract{\endinsert\bigskip}

\def\uspar{\parindent0pt\parskip=\smallskipamount}

\def\supp{\mathop{\rm supp}}

\def\beginproof{{\it Proof: }}
\def\beginproofof#1{{\it Proof of #1: }}
\def\endproof{{\hfill$\square$}\vskip 0pt plus-.1\vsize
\penalty-250\vskip 0pt plus.1\vsize}
\def\title{Asymptotic dimension and uniform embeddings}
\centerline{\bf \title}
\footnote{}{2000 {\it Mathematics Subject Classification: }20F69, 20F65, 20H15; 20E22, 54F45, 51F99.}
\smallskip
\centerline{\'Swiatos\l aw R. Gal%
\footnote{$^\star$}{Partially supported by an
FNRS grant 20-109130 U~00998 and KBN grant 2 P03A 017 25.}}
\centerline{Universit\'e de Neuch\^atel and Wroc\l aw University}
\font\tt=cmtt10
\centerline{\tt http://www.math.uni.wroc.pl/\~{}sgal/papers/asdim.ps}

\beginabstract
We show that the {\sevensl type function} of a space with finite asymptotic
dimension estimates its Hilbert (or any $\ell^p$) compression. The method
allows to obtain the lower bound of the compression of the lamplighter
group $\Bbb Z\wr\Bbb Z$, which has infinite asymptotic dimension.
\endabstract
\uspar

\def\dist{\hbox{\rm dist}}
\def\diam{\hbox{\rm diam}}
\def\asdim{\hbox{\rm asdim}}
\def\U{\frak U}

\section Introduction

The study of embeddings of groups into functional spaces (where
Hilbert space plays a prominent role) was introduced by M.~Gromov in
\cite[Sec.~7.E]{MR1253544}. It appeared to be a fundamental tool
in geometric group theory since G.~Yu proved that a group admitting
a uniform embedding into a Hilbert space satisfies the Novikov Conjecture
on homotopy equivariance of higher signatures \cite{MR1728880} (which
was predicted by M.~Gromov).

In \cite[p.~29]{MR1253544} M.~Gromov also introduced a large scale twin of
the topological covering dimension, the asymptotic dimension.
As noticed by N.~Higson and J.~Roe \cite{MR1739727} spaces of finite asymptotic
dimension have property A of G.~Yu, a thus, in particular, they
uniformly embed into $\ell^p$ spaces. This ideas already
appear in \cite[Sec.~6]{MR1626745}.

In \cite{MR2160829} E.~Guentner and J.~Kaminker initiated the study of rates
of such embeddings. This topic was developed in a number of papers
\cite{ags-metrics,ctv-growth,sv-wreath,tessera-compression,tessera-subexp}.

In this paper we provide an estimate for the compression of such an
embedding (see \ref{T:main}) in terms of the secondary asymptotic invariant
of a space with finite asymptotic dimension, the type function (see
\ref{S:type} for the definitions).

We will also use those estimates to get the lower bound
of the compression of the lamplighter
group $\Bbb Z\wr\Bbb Z$, which has infinite asymptotic dimension.

\subsection Results

The {\it type function} $D_k$ (\ref{S:type}, cf. \cite[p.~29]{MR1253544})
of a metric space $X$ is defined as follows. $D_k(L)$ is the infimum
of those $S$ such that there exist an open cover $\U$ of $X$
by sets of diameter almost $S$, the multiplicity of $\U$ is
at most $k+1$ and, for every $x\in X$, the ball $B_L(x)$ of radius $L$
around $x$ is contained in some set form $\U$. The space $X$ is said
to have {\it asymptotic dimension at most $k$} if
$D_k(L)<\infty$ for all positive $L$.

\thm Theorem. \label{T:main}
Let $(X,d)$ be a space of bounded geometry, asymptotic dimension at
most $k$ and {\rm type function} $D_k$.
Let $u$ be any non-decreasing function such that
$$\int_c^\infty{du(t)^p\over t^p}<\infty.\label{Eq:tessera-cond}$$
Then there exists a Lipschitz map $\theta\colon
X\to L^p([c,\infty),du^p)\otimes\ell^p(X)$ such that
$$d(x,y)\geq D_k(t)\ \Rightarrow\ ||\theta(x)-\theta(y)||_p\geq u(t).$$

In particular, one can take $u(t)=t\cdot(\log t)^{-(1+a)/p}$ for any $a>0$
(see \ref{E:overlog}).
Other form of the condition \ref{Eq:tessera-cond} appears in
\cite{tessera-subexp} as condition ($C_p$).

The proof of \ref{T:main} will be postponed until \xref{proof:main}
after \ref{C:epsdel}. 

If $(X,d)$ is a group with the word metric then the type function
grows at least linearly since the diameter of the ball of radius $L$
is at least $\lfloor L\rfloor$. The same is true for quasigeodesic spaces.

\rmk Examples. The class of spaces with linear type function contains
trees and Loba\v cevski\u{\i} hyperbolic spaces \cite[p.~29]{MR1253544}.

Unfortunately, the most common way to show that a space has linear type
function is to embed it into a product of trees and/or hyperbolic spaces.
M.~Bonk and O.~Schramm \cite{MR1771428}
showed that every Gromov hyperbolic group admits
a quasi-isometric embedding into
a Loba\v cevski\u{\i} hyperbolic space.
Alternatively, S.~Buyalo and V.~Schroeder \cite{MR1707192}
shoved that such a group embeds quasi-isometrically in
a product of a finite number of (locally finite) regular trees.

Also Coxeter groups have linear type function as they can
be embedded quasi-isometrically in a product of (locally finite) regular
trees \cite{MR1925487}.

Another examples are amenable Baumslag-Solitare groups defined by
the presentation $BS^1_m=\langle a,t|tat^{-1}=a^m\rangle$ (where
$m>1$ is any integer). They embed quasi-isometrically
into the product of a $(m+1)$-valent
tree and a Loba\v cevski\u{\i} space \cite{MR1608595}.

Recently, P.~Nowak \cite{nowak-exact} found the first examples of groups
with finite asymptotic dimension (even of asymptotic dimension one)
with nonlinear type function.

\subsection Acknowledgements

The author warmly
thanks Piotr Nowak and Romain Tessera for
inspiring discussions.
The author is also grateful to Pierre de la Harpe for
hospitality during the final work on the manuscript.

\subsection Notation and preliminaries

Let $f,g\colon\Bbb R_+\to\Bbb R_+$ be two weakly monotone functions.
We write $f\succeq g$ if there exist positive constants $C$ and $D$
such that $f(t)\geq Cg(Dt)$ for sufficiently large $t>0$. We write
$f\sim g$ if $f\succeq g$ and $f\preceq g$.

If $f$ is a non-decreasing (resp. non-increasing) positive function, then
by $f^{-1}(t)$ we mean $\inf f^{-1}([t,\infty))$
(resp. $\sup f^{-1}((0,t])$). 

\thm Definition. A metric space is said to have {\rm bounded geometry}
if for every $R\geq0$ there exists $C<\infty$ such that every ball of
radius $R$ contains at most $C$ points.

Most of this work remains true for a wider class of spaces
namely {\it measure-metric spaces with bounded geometry}
(cf. \cite{tessera-subexp}).

By $\ell^p_1X$ we will denote the unit sphere in $\ell^pX$.

\section Uniform embeddings in $\ell^p$-spaces

\subsection Property A and its profile

In this section we introduce a quantitative description of a Property A
of G.~Yu which will be subsequently used in the proof of \ref{T:main}.
Other definition were recently proposed by P.~Nowak
\cite[Def.~3.2]{nowak-exact} and R.~Tessera \cite[Def.~3.1]{tessera-subexp}.

\thm Definition.
Let $\xi\colon X\times X\to\Bbb R$ be a 
kernel. We will write $(x,y)\mapsto \xi_x(y)$. Define
{\parindent=.3in\parskip=0pt
\item{(1)} $S(\xi):=\sup\{d(x,y)\colon \xi_x(y)\neq 0\}$,
\item{(2)} $\varepsilon(\xi;p):=
\sup\left\{{||\xi_x-\xi_y||_p\over d(x,y)}|x\neq y\right\}$.}

\rmk Remark. A map $\xi\colon X\to\ell^p(X)$ is
$\varepsilon(\xi;p)$-Lipschitz,
and $\varepsilon(\xi;p)$ is the best Lipschitz constant.

\thm Definition. For a space $X$ define
$$\epsilon_{X;p}(S)=\inf\{\varepsilon(\xi;p)|\xi\colon X\to \ell^p_1(X),\ 
S(\xi)\leq S\}.$$
\label{D:profile}

In particular, for any $S>0$, it is possible to find a map
$\xi\colon X\to \ell^p_1(X)$ with $S(\xi)\leq S$ and
$\varepsilon(\xi;p)<2\epsilon_{X;p}(S)$.
Notice that $\epsilon_{X;p}$ is a non-increasing function.
We will write $\epsilon_{X;p}=:\epsilon_p$ if it does not lead to
ambiguity.

We do not recall the original definition of the property A due to
G.~Yu \cite[Def.~2.1]{MR1728880}. Instead we give an equivalent formulation
by N.~Higson and J.~Roe \cite[L.~3.2]{MR1739727}.

\thm Definition. A discrete metric space $X$ of bounded geometry
has property A if for any $R>0$ and $\epsilon>0$ there exists
$\xi\colon X\to\ell_1^1X$ with $S(\xi)<\infty$ and such that
$||\xi_x-\xi_y||\leq\epsilon$ provided $d(x,y)\leq R$.

\thm Proposition. If\/ $\inf\limits_S\epsilon_{X;1}(S)=0$
then $X$ has property A.
\label{P:limA}

\beginproof
$\xi\colon X\to\ell_1^pX$ with $S(\xi)<\infty$ and
$\varepsilon(\xi;1)=\epsilon/R$ satisfies
the condition in the above definition.
\endproof

To prove the opposite implication we need a mild assumption that the space
is uniformly discrete (which can be always realized in a
quasi-isometry class of the metric).

\thm Definition. A metric space $(X,d)$ is {\rm uniformly discrete} if
$0$ is an isolated value of the metric.

An example of uniformly discrete space is a vertex set of
a graph (with unit length edges) with an induced metric,
eg. a discrete group with a word metric.

\thm Proposition. Assume that a uniformly discrete space $X$ has property A.
Then $\inf\limits_S\epsilon_{X;1}(S)=0$.

\beginproof
By hypothesis we can find $r>0$ such that $d(x,y)<r$ implies $x=y$.
Assume that $||\xi_x-\xi_y||\leq\epsilon\cdot r$ for $d(x,y)\leq2/\epsilon$.
Then $\xi\colon X\to\ell^1(X)$ is $\epsilon$-Lipschitz.
\endproof

We are interested in asymptotic behavior of $\epsilon_{X;p}$
(how fast it does converge to zero).

\rmk Example.
Let $V$ be a vertex set (with the induced metric) of any simplicial tree.
Then $\epsilon_{V;p}(S)\geq(2S)^{-1/p}$.
Indeed, fix a point  $\omega$ in the boundary of the tree.
Let $\xi_x(y)=S^{-1/p}$ when $y$ is at distance at most $S$
from $x$ in the direction of $\omega$. Then $\xi\colon V\to \ell^p(V)$
is $(2S)^{-1/p}$-Lipschitz \cite{MR1802681}.
\label{E:badtree}

The above estimate is sharp only for $p=1$.
The following example shows the optimal estimate.
Although it is a direct corollary of \ref{C:epsdel}
in the case where the tree is uniformly locally finite (has bounded geometry),
we find it instructive to do the proof by hand in full generality.

\rmk Example.
Let $V$ be a vertex set (with the induced metric) of any simplicial tree.
For any $S\geq 1$ one can construct $\xi\colon V\to \ell^p(V)$
with $S(\xi)\leq S$ and such that $\xi$ is $8S^{-1}$-Lipschitz.
\label{E:goodtree}

\beginproof
Indeed, define $\zeta_x(z)=\max(S+2-|S-2d(x,z)|,0)$ if the geodesic
ray from $x$ towards $\omega$ goes through $z$ and $\zeta_x(z)=0$ otherwise.
Then
$$||\zeta_x||_p>\left(2\int_0^{S/2}(2t)^p\,dt\right)^{1/p}=
\left({S^{p+1}\over p+1}\right)^{1/p}
={S^{1+1/p}\over(p+1)^{1/p}},$$
and if $d(x,y)=1$, then $||\zeta_x-\zeta_y||_p^p=2^p(2\lfloor S/2\rfloor+2)$.
Thus $\xi_x=\zeta_x/||\zeta_x||_p$ satisfies
$$||\xi_x-\xi_y||_p<2\left({2\lfloor S/2\rfloor+2\over S}\right)^{1/p}
(1+p)^{1/p}S^{-1}\leq 8S^{-1}$$
for $d(x,y)=1$.
\endproof

\subsection Quasi-isometry invariance

Let $f\colon X\to Y$ be a map. If $X$ and $Y$ are
equipped with metrics $d_X$ and $d_Y$ we define the compression $\rho_f$
of $f$ as the greatest non-decreasing function such that
$$\rho_f(d_X(x,x'))\leq d_Y(f(x),f(x')).$$

\thm Proposition. \label{P:subspace}
Let $f\colon X\to Y$ be a map. Then
$\epsilon_{Y;p}(\rho_f(S))\leq\epsilon_{X;p}(3S)$.

\beginproof
Chose a map $s\colon Y\to Y$ such that $s(y)\in f(X)\cap B_y(2\dist(y,f(X))$.
Obviously $f(X)\cap B\left(y,2\dist(y,f(X)\right)$ is not empty.
Notice that
$$d(f(x),s(y))\leq d(f(x),y)+d(s(y),y)\leq 3d(f(x),y)$$
while $d(s(y),y)\leq 2\dist(f(X),y)\leq 2d(f(x),y)$.
In particular, $s(y)=y$ if $y\in f(X)$.

For $\xi\colon Y\to\ell^p(Y)$ we define $\sigma\colon X\to\ell^p(X)$
by the formula
$$\sigma_x(z):=\left(\sum_{s(y)=f(z)}|\xi_{f(x)}(y)|^p\right)^{1/p}.$$

We are left to show that
$$\eqalign{||\sigma_x||_p&=||\xi_{f(x)}||_p,\cr
||\sigma_x-\sigma_{x'}||_p&\leq ||\xi_{f(x)}-\xi_{f(x')}||_p,\rlap{ and}\cr
\rho_f(S(\sigma))&\leq3S(\xi).\cr}$$

Indeed,
$$\eqalign{||\sigma_x||_p^p&=\sum_{z\in X}\sum_{s(y)=f(z)}|\xi_{f(x)}(y)|^p
=\sum_{y\in Y}|\xi_{f(x)}(y)|^p=||\xi_{f(x)}||_p^p.\cr
||\sigma_x-\sigma_{x'}||_p^p&=\sum_{z\in X}
\left(\left(\sum_{s(y)=f(z)}|\xi_{f(x)}(y)|^p\right)^{1/p}
-\left(\sum_{s(y)=f(z)}|\xi_{f(x')}(y)|^p\right)^{1/p}\right)^p\cr
\noalign{\hbox{(by triangle inequality in $\ell^p(s^{-1}\circ f(z))$)}}
&\leq\sum_{z\in X}\sum_{s(y)=f(z)}
\left|\xi_{f(x)}(y)-\xi_{f(x')}(y)\right|^p=||\xi_{f(x)}-\xi_{f(x')}||_p^p.\cr
\rho_f(S(\sigma))&=\sup\{\rho_f(d_X(x,z))|\sigma_x(z)\neq0\}\cr
&\leq\sup\{d_Y(f(x),f(z))|\sigma_x(z)\neq0\}\cr
&\leq\sup\{d_Y(f(x),s(y))|\xi_{f(x)}(y)\neq0\}\cr
&\leq\sup\{3d_Y(f(x),y)|\xi_{f(x)}(y)\neq0\}\cr
&\leq\sup\{3d_Y(u,y)|\xi_u(y)\neq0\}=3S(\xi).\cr}$$
\endproof

A consequence of the previous proposition is the quasi-isometry
invariance of $\epsilon_{X;p}$.

\thm Corollary.
Let $X,Y$ be metric spaces. Let $f:X\to Y$ be a quasi-isometry. Then
$$\epsilon_{X;p}\sim\epsilon_{Y;p}.$$

From \ref{P:subspace} we get an immediate
\thm Corollary.
Let $X\subseteq Y$ be a subspace with the induced metric. Then
$$\epsilon_{X;p}\preceq\epsilon_{Y;p}.$$

\subsection Dependence on $p$

For the proofs of the estimates of this section the reader
may consult \cite{MR1727673} or easily adjust the proofs
from the original paper \cite{Mazur}.

Consider a map called the Mazur map $M\colon \ell^q_1X\to\ell^p_1X$
defined by the formula
$$(Mf)(x)=|f(x)|^{q/p-1}f(x).$$

This map is $q/p$-Lipschitz if $p\leq q$. Thus we obtain
\thm Corollary.  Assume that $p\leq q$. Then
$$\epsilon_p\leq q/p\epsilon_q.$$
\label{C:epsLipsch}

\thm Lemma.
Assume that $\epsilon_p(S)\leq \varphi(S)S^{\alpha/p}$, where $\varphi$
does not depend on $p$. Then
$$\epsilon_p(S)\leq {e^\alpha\over p}\varphi(S)\log(S)$$
for $S\geq e^p$.
\label{L:log}

\beginproof
This follows form \ref{C:epsLipsch} by putting $q=\log S$.
\endproof

On the other hand if $p\geq q$ the Mazur map is only H\"older
with exponent $p/q$.

Let us make a reasonable assumption on the metric space.
We assume that there exists $M$ such that for $d$
sufficiently large there are $x,y\in X$ such that $d<d(x,y)<Md$.

Under the above assumption $\epsilon_p(S)\leq (MS)^{-1}$.
Indeed if $d(x,y)>2S$ and $S(\xi)\leq S$
then $\xi_x$ and $\xi_y$ have disjoint supports and $||\xi_x-\xi_y||=2$.
Thus if $2MS>d(x,y)>2S$ then $||\xi_x-\xi_y||/d(x,y)>(MS)^{-1}$.

Thus even if $\epsilon_1(S)\preceq S^{-1}$ (i.e. $\epsilon_1$ has the fastest
possible decay) the estimate on $\epsilon_p$
coming from the Mazur map ($\epsilon_p(S)\preceq S^{-1/p}$)
is usually far from sharp (cf. \refs{E:badtree} and \refn{E:goodtree}).

\subsection $\ell^p$ compression

\thm Definition. \cite{MR1253544}
Let $X,Y$ be metric spaces. A map $f\colon X\to Y$ is
a {\rm coarse embedding (uniform map)} if there exist non-decreasing functions
$\rho_-,\,\rho_+\colon [0,\infty)\to[0,\infty)$ satisfying
{\parindent=.3in\parskip=0pt
\item{(1)}$\rho_-(d_X(x,y))\le d_Y(f(x),f(y))\le \rho_+(d_X(x,y))$
for all $x,y\in X$,
\item{(2)}$\lim_{t\to\infty}\rho_-(t)=+\infty$.}

The smallest function $\rho_+$ one can choose is called the
{\sl dilation} of $f$.

If the space $X$ is quasi-geodesic
(as for example a Cayley graph of a group)
then one can take $\rho_+$ to be an affine function.
The main quantitative interest is how
big can one choose $\rho_-$.

Hilbert space compression rate of a metric space $X$ was
introduced by Guentner and Kaminker \cite{MR2160829} and it is the
supremum of $\alpha$ such that $X$ admits a coarse embeddings into
the Hilbert space with lower bound $\rho_-(t)\succeq t^{\alpha}$.

Below we consider embeddings in an $\ell^p$-space for any $p\geq 1$.

\thm Proposition. Let $X$ be a metric space.
Let $f\colon [c,\infty)\to\Bbb R$ be a non-decreasing left-continuous function.
Assume that a measurable field of maps
$$[c,\infty)\times X\ni (S,x)\mapsto\xi^S_x\in\ell^p_1X$$
satisfies $S(\xi^S)\leq S$ and
$$\int_c^\infty\varepsilon(\xi^S,p)^p df(S)^p=:C<\infty.$$
Then the map $\theta(x):=\xi(x)-\xi(x_0)$ (where $x_0$ is an arbitrary
reference point) is a $C^{1/p}$-Lipschitz map
$X\to L^p([c,\infty),df^p)\otimes\ell^p(X)$
which compression $\rho_-$ is asymptoticly bounded by $f$, more precisely
$$\rho_-(d)\geq 2f(d/2)-2f(c).$$
\label{P:compression}

\rmk Remark. One can always find $f$ piecewise constant with the same
asymptotic behavior.

\beginproof
{\it Step 1 ($\theta$ is C-Lipschitz)} By the definition of the norm in
$L^p([c,\infty),df^p)\otimes\ell^p(X)$ we have
$$||\theta(x)-\theta(y)||_p^p=\int_c^\infty||\xi^S_x-\xi^S_y||_p^pdf(S)^p
\leq\left(\int_c^\infty\varepsilon(\xi^S,p)^p df(S)^p\right)d(x,y)^p=
C\,d(x,y)^p.$$

{\it Step 2 (bound on the compression).} Assume that $d(x,y)>d$.
Since $\supp(\xi^S_x)\subset B_S(x)$ w have that $||\xi^S_x-\xi^S_y||_p=2$
if $2S\leq d$. What follows,
$$\eqalign{||\theta(x)-\theta(y)||_p^p&\geq \int_c^{d/2} 2^p df(S)^p\cr
&=2^p(f(d/2)^p-f(c)^p)\geq 2^p(f(d/2)-f(c))^p.\cr}$$

The claim follows by the continuity of $f$ from the left.
\endproof

\thm Corollary. 
Let $X$ be a metric space.
Let $f\colon [c,\infty)\to R$ be a non-decreasing function such that
$$\int_c^\infty\epsilon_p(S)^p df(S)^p<\infty.$$
then there exist an uniform embedding $\theta\colon X\to
L^p([c,\infty),df^p)\otimes\ell^p(X)$ such that the compression $\rho_-$
satisfies $\rho_-\succeq f$.
\label{C:compression}

\beginproof
The assumption in \ref{P:compression} that $f$ is left-continuous
was made only to get a precise bound on the compression. Replacing
$f$ by $g(t):=\lim_{s\searrow t}f(s)$ makes the function $g$ left-continuous
and does not change the value of the integral. On the other hand $f\sim g$
with constants arbitrary close to $1$.

Taking sufficiently small subdivision of $[c,\infty)$ one constructs
a piecewise constant (thus measurable) field $\xi$ with the property
$\varepsilon(\xi^S;p)<2\epsilon_p(S)$ and $S(\xi^S)\leq S$.
Thus the claim.
\endproof

\rmk Example.
Let $u$ satisfy the condition
$$\int_c^\infty{du(t)^p\over t^p}<\infty.$$
Then $f(t):=u(1/\epsilon_p(t))$ satisfies the assumption of
\ref{P:compression}. An example of such $u$ is $u(t)=t\log(t)^{-(1+a)/p}$
for any $a>0$. Indeed,
$$\eqalign{\int_c^\infty\epsilon_p(S)^p
d{1\over \epsilon(S)^p|\log(\epsilon_p(t))|^{1+a}}&\cr
&[t:=|\log(\epsilon_p(S))|]\cr
=\int_{|\log\epsilon_p(c)|}^\infty e^{-pt}d{e^{pt}\over t^{1+a}}
&=\int_{|\log\epsilon_p(c)|}^\infty d\,t^{-1-a}+pt^{-1-a}d\,t\cr
&\left.t^{-1-a}-{p\over a}t^{-a}\right|_{|\log\epsilon_p(c)|}^\infty
<\infty.\cr}$$
\label{E:overlog}

\subsection Uniform embeddings of trees

\ref{C:compression} and \ref{E:goodtree} provide an embedding of a tree
into $\ell^p$ space with compression bigger than given non-decreasing
function $u$ satisfying $\int t^{-p}du(t)^p<\infty$.

Independently, Tessera improved \cite[Thm.~7.3]{tessera-compression}
the original construction of Guentner and Kaminker \cite[Prop.~4.2]{MR2160829}
to obtain an embeddings of a trees with such an asymptotic.

On the other hand Tessera \cite[Cor.~6.3]{tessera-compression}
showed, that, for $2\leq p<\infty$ and any tree $T$ with no vertices
of valence $1$ or $2$,
the compression $\rho$ of any Lipschitz map $T\to\ell^p$ satisfies
$\int t^{-p}d\rho(t)^p<\infty$.

The difference between our construction and the construction of
Guentner, Kaminker and Tessera is the following. The former is
{\sl a cocycle} with respect to the action of the (amenable) stabilizer
of the point in the boundary, when the latter is a cocycle with respect
to the the action of the (compact) stabilizer of the vertex in the tree.

By being a cocycle we mean the following property. Let a group $\Gamma$
act on a space $X$. The action induces a representation on the space
$W=\ell^p(X)\otimes L^p([c,\infty))$. The map $\theta\colon X\to W$ is called
$\Gamma$-{\sl cocycle} if
$\Gamma\ni\gamma\mapsto\vartheta(\gamma):=
\theta(\gamma x)-\gamma\theta(x)\in W$
is independent on $x\in X$. Moreover then $\vartheta$ satisfies
the cocycle relation
$$\vartheta(\gamma\gamma')=\gamma\vartheta(\gamma')+\vartheta(\gamma),$$
or, in other words, $w\mapsto \gamma w+\vartheta(\gamma)$ is an affine
isometric action of $\Gamma$ on $W$.

\thm Proposition. Let a group $\Gamma$ acts on a space $X$.
If $\xi^S\colon X\to\ell^p_1X$ are $\Gamma$-equivariant
then the map $\theta$ constructed in \ref{P:compression} is a $\Gamma$-cocycle.
Moreover, $\vartheta(g)=\theta(gx_0)$ where $x_0$ is the chosen reference
point.

On the other hand, assume that $\theta$ is an cocycle on the vertex set of
a tree with values in a Hilbert space equivariant with respect to some group
$\Gamma$, subgroup of the full automorphism group of the tree,
with the compression $\rho_\theta$ satisfying
$\rho_\theta(t)/\sqrt t\to\infty$. Then the group $\Gamma$ is necessarily
{\sl amenable} by \cite[Thm.~4.1]{ctv-growth} or by a small modification of
\cite[Thm.~5.3]{MR2160829}.

\section Spaces with finite asymptotic dimension

\subsection Type function
\label{S:type}

\thm Definition.
Let $\U$ be a cover of $X$. Define
{\parindent=.3in\parskip=0pt
\item{(1)} the Lebesgue number at $x\in X$
$$L(\U,x)=\sup_{U\in\U}\{r\colon B_r(x)\subset U\},$$
\item{(2)} the Lebesgue number of $\U$
$$L(\U)=\inf_{x\in X}L(\U,x),$$
\item{(3)} the multiplicity of $\U$ at $x\in X$
$$m(\U,x)=\#\{U\in \U|x\in U\},$$
\item{(4)} the multiplicity of $\U$
$$m(\U)=\max_{x\in X}m(\U,x),$$
\item{(5)} the mesh of $\U$
$$S(\U)=\sup\{\diam(U)|U\in\U\}.$$}

\thm Definition. \cite{db}
We say that a metric space has asymptotic dimension less than
$m\in\Bbb N$, denoted $$\asdim X<m,$$ if for every $L<\infty$ there exist
a number $S<\infty$ and a cover $\U$ with mesh at most $S$,
with multiplicity at most $m$ and Lebesgue number at least $L$.

\thm Definition. 
Let $X$ be a metric space with finite asymptotic dimension. Define
the type function $D_{m-1}\colon\Bbb R^+\to\Bbb R^+\cup\{\infty\}$
in the following way:
$D_{m-1}(L)$ is the infimum of those $S>0$ for which $X$ can be covered
by a family of sets with mesh at most $S$,
multiplicity at most $m$ and Lebesgue number at least $L$.


\thm Definition. We also define
$\delta_p(S)=\sup_\U\left\{{L(\U)\over m(\U)^{2/p}}\right\}$
where $\U$ runs over a set of covers of $X$ with mesh at most $S$.
\label{D:type}

\rmk Remark. Of course, $D_{m-1}(m^{2/p}\delta_p(S))\geq S$
for any $m$ and $p$.
\label{R:delta}


\thm Remark. The type function was originally (\cite[p.~29]{MR1253544})
defined in a different way. In the rest of this section we will compare
the two definitions.

\thm Proposition.
A metric space $X$ has asymptotic dimension at most $k$
if for every $L$ there exists a cover $\U=\bigcup_{i=0}^k\U_k$ with finite
mesh and the property, that, for any $0\leq i\leq k$, any two different
sets form $\U_i$ are $L$ disjoint.
\label{P:orig}

\beginproof
This is a part of \cite[Thm.~1]{db}.
\endproof

The original definition of a type function, which we will call
$\widetilde D_k$, is as follows: $\widetilde D_k(L)$ is the infimum
of the mesh of the covers as in the \ref{P:orig}.

Note that $L/2$ thickening of the cover as in \ref{P:orig} is a cover
with mesh at most $\widetilde D_k(L)+L$ multiplicity at most $k$
and Lebesgue number at least $L/2$. Thus
$$D_k(L/2)\leq\widetilde D_k(L)+L.$$
We leave it as an easy exercise to show that if the space is the vertex set
of a graph (with unit length edges) with the induced metric, then
$\widetilde D_k(L)\geq L/k-1$. More generally, the inequality
$\widetilde D_k(L)\succeq L$ holds for quasi geodesic spaces.
For such spaces $D_k\preceq\widetilde D_k$ and
\ref{T:main} remains true if one replaces $D_k$
by $\widetilde D_k$.

\subsection Asymptotic dimension and Property A

In this section we will exhibit a link between $\epsilon_p$ and asymptotic
dimension of $X$, in particular we will show that a large class of spaces,
namely the spaces with finite asymptotic dimension of {\sl linear type},
satisfy $\epsilon_{X;p}(S)\sim S^{-1}$ for all $p\geq 1$.

\thm Theorem.
Let $\U$ be a cover of\/ $X$ with finite multiplicity. Then there exist a map
$\xi\colon X\to\ell^p_1X$ which is 
$2{\root p\of{2m(\U)^2}\over L(\U)}$-Lipschitz
and satisfies $S(\xi)=S(\U)$.
\label{T:coverLipsh}

\beginproof
Let $\U$ be a cover of $X$. For any $U\in\U$ define
$$\eqalign{\psi_U(x)&=\dist(x,X-U),\cr
\psi(x)&=\left(\sum_{U\in\U}\psi_U(x)^p\right)^{1/p}.\cr}$$

We have
$$\eqalign{\psi(x)&\geq L(\U,x)\geq L(\U),\hbox{\rm\ and}\cr
\psi_U&\hbox{\rm\ is $1$-Lipshitz, i.e. }|\psi_U(x)-\psi_U(y)|\leq d(x,y).\cr}$$

{\sl Claim:}
$$\left|{\psi_U(x)\over\psi(x)}-{\psi_U(y)\over\psi(y)}\right|
\leq 2{(m(\U)+1/2)^{1/p}\over L(\U)}d(x,y).$$
Indeed,
$$\eqalign{\left|{\psi_U(x)\over\psi(x)}-{\psi_U(y)\over\psi(y)}\right|^p&=
\left|{\psi_U(x)-\psi_U(y)\over\psi(x)}
-{\psi_U(y)\over\psi(y)}\cdot{\psi(x)-\psi(y)\over\psi(x)}
\right|^p\cr
\noalign{\hbox{(by the inequality between the arithmetic and $p$-mean)}}
&\leq2^{p-1}\left({|\psi_U(x)-\psi_U(y)|^p\over\psi(x)^p}
+{\psi_U(y)^p\over\psi(y)^p}\cdot{|\psi(x)-\psi(y)|^p\over\psi(x)^p}\right)\cr
\noalign{\hbox{(by the triangle inequality in $\ell^p(\U)$ and
$\psi_U(y)\leq\psi(y)$)}}
&\leq{2^{p-1}\over\psi(x)^p}\left(|\psi_U(x)-\psi_U(y)|^p+\sum_{V\in\U}
|\psi_V(x)-\psi_V(y)|^p\right).\cr
&\leq{2^{p-1}\over L(\U)^p}(1+2m(\U))d(x,y)^p.\cr}$$

Let $\phi_U(x)=\psi_U(x)/\psi(x)$ and $\chi_U=\psi_U/||\psi_U||_p$.
Define $$\xi_x(z):=
\left(\sum_{U\in\U}\phi_U(x)^p\chi_U(z)^p\right)^{1/p}.$$
Note that $\chi_U$ and $\xi_x$ are well defined.
The support of $\psi_U$ is contained in a ball of radius $S$, thus the
norm reduces to the finite sum. Also the sum in the definition of $\xi_x(z)$
runs over those $U$ which contain both $x$ and $z$, thus at most $m(\U)$
of them. In particular, $\xi_x(z)\neq 0$ if and only if
$x$ and $z$ belong simultaneously to at least one $U\in\U$.
Therefore $S(\xi)=S(\U)$.

Moreover,
$||\xi_x||_p^p=\sum_{z\in X}\sum_{U\in\U}
\phi_U(x)^p\chi_U(z)^p=
\sum_{U\in\U}\phi_U(x)^p\left(\sum_{z\in U}\chi_U(z)^p\right)=1$.

We are left to check the Lipschitz condition. If $d(x,y)\geq L(\U)$,
the condition is trivial as $||\xi_x-\xi_y||\leq||\xi_x||+||\xi_y||=2$
for all $x$ and $y$. Therefore assume that $d(x,y)<L(\U)$, and, in particular,
there is a set in $\U$ containing both $x$ and $y$. Thus, there
are at most $2m(\U)-1$ sets containing any of $x$ or $z$.
As previously we check:
$$\eqalign{||\xi_x-\xi_y||_p^p&=\sum_{z\in X}\left|\left(\sum_{U\in\U}
\left(\phi_U(x)\chi_U(z)\right)^p\right)^{1/p}-\left(\sum_{U\in\U}
\left(\phi_U(y)\chi_U(z)\right)^p\right)^{1/p}\right|^p\cr
\noalign{\hbox{(by triangle inequality in $\ell^p(\U)$)}}
&\leq\sum_{z\in X}\sum_{U\in\U}\left|\phi_U(x)\chi_U(z)-
\phi_U(y)\chi_U(z)\right|^p\cr
&=\sum_{U\in\U}\left|\phi_U(x)-\phi_U(y)\right|^p
\left(\sum_{z\in U}\chi_U(z)^p\right)\cr
\noalign{\hbox{(by the Claim)}}
&\leq\sum_{x\hbox{\sevenrm\ or }y\in U\in\U}
2^{-1}(2m(\U)+1)\left({2d(x,y)\over L(\U)}\right)^p\cr
&=2^{-1}(2m(\U)-1)(2m(\U)+1)\left({2d(x,y)\over L(\U)}\right)^p
\leq2m(\U)^2\left({2d(x,y)\over L(\U)}\right)^p.\cr}$$
and, in particular, $\xi\colon X\to \ell_1^p(X)$ is
$2(2m(\U)^2)^{1/p}L(\U)^{-1}$-Lipschitz.
\endproof

The proof of \ref{T:coverLipsh}
depends on an informal argument of Higson and Roe \cite{MR1739727}
and more precisely on a computation from the proof of Theorem~1
form \cite{db} (case $p=1$ of \ref{T:coverLipsh}).

\thm Corollary.
Let $X$ be a metric space with finite asymptotic dimension and
$\delta$-function $\delta_p$ (see \ref{D:type}). Then
$$\epsilon_p\leq {2^{1+1/p}\over\delta_p}\leq{4\over\delta_p}$$
for all $1\leq p<\infty$.
\label{C:epsdel}

\beginproofof{\ref{T:main}} \xrdef{proof:main}
Substitute $f:=u\circ D_k^{-1}$ in \ref{C:compression}. The claim
follows from inequalities
$$4/\epsilon_p\geq\delta_p\geq k^{-2/p}D_k^{-1},$$
due to \ref{R:delta} and \ref{C:epsdel}.
\endproof

A.~N.~Dranishnikov \cite{dran-poly} defined groups with a polynomial
dimension growth. In term of the function $\delta_p$ it is defined as follows

\thm Definition. 
A space has {\rm polynomial dimension growth of degree less than} $k$
if\/ $\lim\limits_{S\to\infty}\delta_{2k}(S)=\infty$.

A straightforward corollary from \ref{C:epsdel} and \ref{P:limA} is

\thm Corollary. \cite[Thm.~3.3]{dran-poly} A space with a polynomial
dimension growth has property~A.

\rmk Example. A space $X$ is said to have polynomial growth if there exist $C$
and $n$ such that every ball of radius $R$ contains at most $C\cdot R^n$
elements. Take a cover of $X$ by all balls of radius $S/2$. Its mesh is
at most $S$, Lebesgue number is $S/2$ and multiplicity is $C\cdot (S/2)^n$.
Thus $\delta_p(S)\geq {2^{2/p-1}\over C^{2/p}}S^{1-2n/p}$. In particular, 
by \ref{C:epsdel} $\epsilon_p(S)\leq4C^2S^{-1+2n/p}$, and by
\ref{L:log}, $\epsilon_p(S)\leq(4e^{2n}C^2p^{-1})\log(S)/S$
for $S\geq e^p$, and the $\ell^p$ compression rate for such a space equals one.

\section Application to spaces with infinite asymptotic dimension

\subsection Preliminaries

Let $\Bbb A^{n-1}=\{x=(x_1,\ldots,x_n)\in\Bbb R^n|\sum_{i=1}^nx_i=0\}$.
Equip $\Bbb A^{n-1}$ with a norm $||x||:=\sum_{i=1}^n|x_i|$. 
For any $I\subset\{1,\ldots,n\}$
define a functional $$\phi_I(x):={\sum_{i\in I}x_i\over\#I}-
{\sum_{i\not\in I}x_i\over\#I^c}={n\sum_{i\in I}x_i\over\#I\cdot\#I^c}.$$
Where $I^c$ denotes the complement of $I$. Note that $\phi_I=-\phi_{I^c}$.

\rmk Remark. For any $I$ we have $||x||>\phi_I(x)$.
\label{R:norm}

Let $$V:=\left\{x\in\Bbb A^{n-1}|
\phi_I(x)\leq {1\over2}\hbox{ for all }I\right\}.$$

Let $$[i]_j:=\cases{i-n&if $j\leq i$,\cr i&if $j>i$.\cr}$$
Define $\Lambda:=\{\Bbb A^{n-1}\cap\Bbb Z^n\}$.
Let $\U_i:=\{[i]+x+V|x\in\Lambda\}$.

\thm Lemma. $\U:=\bigcup_{i=0}^{n-1}\U_i$ is a (closed)
cover of $\Bbb A^{n-1}$.

\beginproof
Put for the moment auxiary $\ell^2$-norm,
defined by $||x||_2=\sqrt{\sum_{i=1}^n|x_i|^2}$, on $\Bbb A^{n-1}$.

Below we will just check that $V$ consists of the points that are closer
or at the same distance to $0$ than any vector $[i]$ or its image by the
action of the permutation group $S_n$, which is enough to conclude the claim.

Indeed, the set of points which are closer or at the same distance to $0$
than to $[i]$ is defined by the inequlity
$$\langle x|[i]\rangle\leq1/2||[i]||_2^2,\label{Eq:Voronoi}$$
where $\langle\cdot|\cdot\rangle$ is the scalar product associated
to $||\cdot||_2$.

Straightforward computation shows, that \ref{Eq:Voronoi}
is equivalent to $\phi_{\{1,\ldots,i\}}(x)\leq1/2$.

One can prove a stronger statement, namely that
$V$ is the Voronoi cell of the lattice $\bigcup_{i=0}^{n-1}\Lambda+[i]$
\cite[Ch.~6.6 and 21.3.B]{MR1662447}.
\endproof

\thm Lemma. Each $\U_i$ consists of $1/(n-1)$-disjoint sets.

\beginproof
Let $x$ and $x'$ be two different points in $[i]+\Lambda$.
Chose $i_0$ such that $x_{i_0}>x'_{i_0}$. Let $I=\{i_0\}$.
Then
$$\phi_I(x-x')\geq{n(x_{i_0}-x'_{i_0})\over 1\cdot(n-1)}\geq{n\over(n-1)}.$$
Let $z\in x+V$ and $z'\in x'+V$. Then
$$||z-z'||\geq\phi_I(z-z')\geq\phi_I(x-x')-\phi_I(z-x)
+\phi_I(z'-x')\geq{n\over n-1}-{1\over2}-{1\over2}={1\over n-1}$$
by \ref{R:norm}.
\endproof

\rmk Question. What is the best estimate for disjointness?

The the set of extremal points (the vertices) of $V$ is the orbit
of $\sigma=(2n)^{-1}(1-n,3-n,\ldots,n-1)$ under the permutations
of coordinates \cite[Ch.~6.6 and 21.3.B]{MR1662447}.
From now on assume that $n=2k$ is even.
Then $||\sigma||=(4k)^{-1}\cdot 2\sum_{i=1}^k(2i-1)=k/2$.

Thus we have proved

\thm Lemma. The $1\over2(2k-1)$ open thickening $\U'$ of $\U$ satisfies
$L(\U')\geq{1\over2(2k-1)}$, $S(\U')=k+1/(2k-1)$ and $m(\U'_k)=2k$.

\thm Corollary. $D_{2k;\Bbb Z^k}(L)\leq (2k^2-2k+1)L$.
\label{C:ztok}

\beginproof
Observe that the map $\iota\colon\Bbb Z^k\ni(\ldots,z_m,\ldots)\mapsto
(1/2)(\ldots,z_m,-z_m,\ldots)\in \Bbb A^{2k-1}$ is an isometry.
Thus the induced cover $\iota^*\U'$ satisfies the inequality of the
claim for $L={1\over 2(2k-1)}$.

Since we can precompose $\iota$ with a homotety of $\Bbb A^{n-1}$
the claim is true for arbitrary $L$.

\endproof

\rmk Question. What is the rate of $D_{k;\Bbb Z^k}$?
Can one replace a square in $D_{2k;\Bbb Z^k}(L)\leq 2k^2$
by some lower power? What about the estimates for $D_{m^\alpha}$?
\label{Q:ztom}

One can adjust the results of the next section to prove that
if $D_{cm^\alpha;\Bbb Z^m}\leq Cm^\beta$ then
$\epsilon_{p;\Bbb Z\wr\Bbb Z}\preceq\log(S)/{\root1+\beta\of S}$ and,
what follows, the compression rate of $\Bbb Z\wr\Bbb Z$ is at least
${1\over1+\beta}$. On the other hand
G.~Arzhantseva, V.~Guba and M.~Sapir \cite[Thm.~1.8]{ags-metrics}
showed that the Hilbert compression rate of $\Bbb Z\wr\Bbb Z$
is at most $3/4$. 
It follows that $D_{cm^\alpha}(L)\leq m^\beta\cdot L$
is imposible with $\beta<1/3$. Thus the answer to \ref{Q:ztom}
is nontrivial.

\subsection The Lamplighter group $\Bbb Z\wr\Bbb Z$

The lamplighter group is the (restricted) wreath product $Z\wr Z$, where
$H\wr G$ is defined as a semidirect product $\bigoplus_G H\rtimes G$,
where $G$ acts on $\bigoplus_G H$ permuting the factors.
In other word
$$H\wr G=\langle H,G|[a,[b,g]]=[a,b],\ a,b\in H,\ g\in G-\{e\}\rangle.$$

\thm Proposition \cite[Prop 4.2]{dran-poly}. Let $K$ be a normal subgroup
of $G$. Let $G$ be equipped with a (left invariant) metric. Put on $K$
the restricted metric, and induce the metric on $H=G/K$. Let $k$ and $h$ be
two integers. Then
$$D_{kh-1;G}(L)\leq {4\over 3}D_{k-1;K}(6\, D_{h-1;H}(L)).$$

In our case $G=\Bbb Z\wr\Bbb Z$ and $H=\Bbb Z$. Thus
$D_{2k-1;G}(L)\leq {4/3} D_{k-1;K}(24 L)$,
where $K=\bigoplus_{\Bbb Z}\Bbb Z$ with the restricted metric $d_K$.
This translates to the following statement. Given a cover $\frak V$ of $K$
one constructs a cover $\frak W$ of $\Bbb Z\wr\Bbb Z$ such that
$L({\frak W})\geq L({\frak V})/24$,
$S({\frak W})\leq 4/3S({\frak V})$ and
$m({\frak W})\leq 2m({\frak U})$.

Following Dranishnikov \cite{dran-poly} we decompose
$K:=\bigoplus_{\Bbb Z}\Bbb Z=K_m\oplus K^m$, where
$K_m=\bigoplus_{\{-m+1,-m+2,\ldots,m-1\}}\Bbb Z$ and
$K^m=\bigoplus_{\{\ldots,-m-1,-m,m,m+1,\ldots\}}\Bbb Z$.

\thm Lemma. Let $K/K_m\ni\alpha\to g_\alpha\in K$ be a section and let
$\frak V$ be a cover of $K_m$ then
$\widetilde{\frak V}:=\bigcup_\alpha g_\alpha\frak V$
is a cover of $K$ with the same mesh and multiplicity. Moreover
$L(\widetilde{\frak V})=\min\{L({\frak V}),2m+1\}$.

\beginproof
Observe that any two cosets of $K_m$ are at least $2m+1$ apart.
\endproof

\thm Proposition. Let $f\colon (G,d)\to (G',d')$ be a map, and $\U$ be
a covering of $G'$ then
$$\rho_+(L(f^*\U))\geq L(\U),\quad \rho_-(S(f^*\U))\leq S(\U),$$
where $f^*\U=\{f^{-1}U|U\in\U\}$ is the induced cover of $G$,
and $\rho_\pm$ are the compression and dilation of $f$.
\label{P:induced}

\thm Lemma. The natural homomorphism $j:K_m\to\Bbb Z^{2m-1}$ satisfies
$$d_K(x,y)-4(m-1)\leq d_1(f(x),f(y))\leq d_K(x,y).$$

\beginproof
Recall (\cite{MR2143495}) that if $f\colon\Bbb Z\to\Bbb Z$ is an element of
$K$ (i.e. if $f$ has finite support) then the length of $f$ (with respect to
the metric restricted from $\Bbb Z\wr \Bbb Z$) equals to
$$2\max(\{0\}\cup\{k|f(k)\neq 0\})+2\max(\{0\}\cup\{-k|f(k)\neq 0\})
+\sum_{k\in\Bbb Z}|f(k)|.$$
\endproof

\thm Corollary. There exist a cover ${\frak V}$ of $K_m$ satisfying
$L({\frak V})\geq 2m$, $S({\frak V})\leq 16m^3$ and
$m({\frak V})\leq 4m$.

\beginproof
By \ref{C:ztok} and the previous Lemma we may construct a cover
of $K_m$ satisfying the bounds. The bound on mesh
follows from
$$S({\frak V})\leq D_{4m;\Bbb Z^{2m}}(2m)+4(m-1)\leq
(2(2m)^2-2(2m)+1)\cdot (2m)+4(m-1)\leq 16m^3.$$
\endproof

Let $m=12L$. By the result of Dranishnikov, we construct a cover $\frak W$
of $\Bbb Z\wr\Bbb Z$ with the properties
$L({\frak W})\geq{1\over 24}\cdot 2\cdot 12 L=L$,
$S({\frak W})\leq 4/3\cdot16(12L)^3=36864L^3$ and
$m({\frak W})\leq 2\cdot4\cdot12L=96L$.

\thm Corollary. $\delta_{p;\Bbb Z\wr\Bbb Z}(S)\geq
(96)^{-2/p}\left\lfloor{\root 3\of S/36864}\right\rfloor^{1-2/p}\geq
(1/C)S^{1/3-2/(3p)}$.

\thm Corollary. $\epsilon_{p;\Bbb Z\wr\Bbb Z}(S)\preceq \log(S)/{\root3\of S}$
for any $1\leq p<\infty$.

\beginproof
This follows from \ref{C:epsdel} and \ref{L:log}
\endproof

\thm Corollary. The $\ell^p$ compression rate of Guentner and Kaminker of
$\Bbb Z\wr\Bbb Z$ is at least $1/3$ for all $1\leq p<\infty$. 
\label{C:wr-compr}

\rmk Remark. This result is not sharp.
For $p=2$ (Hilbert compression) by other techniques
G.~Arzhantseva, V.~Guba and M.~Sapir \cite[Thm.~1.8]{ags-metrics} and
independently Y.~Stalder and A.~Valette
\cite[Cor.~4.5a]{sv-wreath} showed that the compression
rate is at least $1/2$. Recently R.~Tessera
showed that the Hilbert compression rate is at least $2/3$
\cite[Cor.~15]{tessera-compression}.

On the other hand P.~Nowak showed \cite[Cor.~4.4]{nowak-exact}
that $\epsilon_{1;\Bbb Z\wr\Bbb Z}\sim S^{-1}$,
thus the $\ell^1$ compression rate of $\Bbb Z\wr\Bbb Z$ is one.
By the H\"older property of the Mazur map the $\ell^p$-compression rate
is at least $1/p$, which gives better estimate that \ref{C:wr-compr}
for $1\leq p<3$.

\bigskip
{\bf References}
\medskip
\bibliography{../../bib/asdim}
\bibliographystyle{alpha}

\bye